\documentclass[12pt,a4paper]{article}

\usepackage{amsthm,amsmath,natbib,amssymb,amsfonts,bm}
\usepackage{epsfig}
\usepackage{placeins}
\usepackage{epsf}
\RequirePackage[dvips]{hyperref}

\def\hang{\hangindent\parindent}
\def\rf{\par\noindent\hang}

\newtheorem*{theorem*}{Theorem}

\theoremstyle{definition}

\theoremstyle{remark}

\begin{document}

\baselineskip=21pt

\begin{center}
 {\bf \Large Note on a paradox in decision-theoretic interval estimation}
\end{center}

\medskip

\begin{center}
{\bf \large Paul Kabaila$^{\textstyle{^*}}$}
\end{center}

\medskip

\noindent{\it $\textstyle{^*}$Department of Mathematics and
Statistics, La Trobe University, Victoria 3086, \newline Australia}

\bigskip
\bigskip

\noindent{\bf Abstract}
\medskip

\noindent Confidence intervals are assessed according to two criteria, namely expected
length and coverage probability. In an attempt to apply the decision-theoretic method to
finding a good confidence interval, a loss function that is a linear combination of the interval
length and the indicator function that the interval includes the parameter of interest has been
proposed. We consider the particular case that the parameter of interest is the normal mean,
when the variance is unknown. Casella, Hwang and Robert, {\sl Statistica Sinica}, 1993,
have shown that this loss function, combined with the standard noninformative prior,
leads to a generalized Bayes rule that is
a confidence interval for this parameter which has ``paradoxical behaviour''. We show
that a simple modification of this loss function, combined with the same prior, leads to
a generalized Bayes rule that is the usual confidence interval i.e.
the ``paradoxical behaviour'' is removed.

\bigskip
\bigskip


\vbox{\vskip 4.5cm}


\noindent $^*$ Corresponding author. Address: Department of
Mathematics and Statistics, La Trobe University, Victoria 3086,
Australia; Tel.: +61-3-9479-2594; fax: +61-3-9479-2466. \newline
{\it E-mail address:} P.Kabaila@latrobe.edu.au.

\newpage

\noindent {\bf 1. Introduction}

\medskip

Suppose that the random vector $X$ has pmf or pdf $f(x|\theta)$, where $x \in {\cal X}$
and $\theta \in \Omega_{\theta}$. Also suppose that either (a) $\psi = \theta$ is the parameter
of interest or (b) $\theta^T = (\psi^T, \tau^T)$ with $\psi$ the parameter of interest and
$\Omega_{\theta} = \Omega_{\psi} \times \Omega_{\tau}$. The decision-theoretic
approach to finding a good point estimator $\delta(X)$ of $\psi$ may be described as follows.
Define the loss function $L(\theta, d)$ for the value $d$ of the estimate of $\psi$, when
the true parameter value is $\theta$. Then define the risk function
$R(\theta, \delta) = E_{\theta} \big(L(\theta, \delta(X)) \big )$,
where $E_{\theta}$ denotes the expectation according to the pmf or pdf $f(x|\theta)$ of $X$.
Choose a prior pdf
$\pi$ (possibly improper) such that minimizing
the posterior expected loss,
with respect to $\delta(x)$ for each $x \in {\cal X}$, yields a good (generalized) Bayes rule estimator. Conditions
for admissibility and for minimaxity of this estimator are well-known (see e.g. Berger, 1985,
Lehmann and Casella, 1998 and Robert, 1994).

Finding a good set estimator $C(X)$ of $\psi$ is much more difficult than finding a good point
estimator of $\psi$. This is because a confidence set $C(X)$ is assessed according two criteria,
namely expected volume and coverage probability. We now have two loss functions and the decision-theoretic
approach does not apply directly. An attempt to apply the decision-theoretic approach is to
define the following loss function, which is a linear combination of the interval
length and the indicator function that the interval includes $\psi$:
\begin{equation}
\label{comp_loss_fn}
L(\theta, C) = \text{vol}(C) - k \, {\cal I}(\psi \in C),
\end{equation}
where $k > 0$ and
\begin{equation*}
{\cal I}({\cal A}) =
\begin{cases}
1 &\text{if } {\cal A} \ \ \text{is true} \\
0 &\text{if } {\cal A} \ \ \text{is false}
\end{cases}
\end{equation*}
for any statement ${\cal A}$. This leads to the risk function
\begin{equation*}
R(\theta, C) = E_{\theta} \big(L(\theta, C(X)) \big)
= E_{\theta}(\text{vol}(C(X))) - k P_{\theta}(\psi \in C(X)),
\end{equation*}
where $P_{\theta}$ denotes the probability according to the pmf or pdf $f(x|\theta)$ of $X$.
One then seeks $k$ and prior pdf $\pi$ such that minimizing
the posterior expected loss, with respect to $C(x)$ for each $x \in {\cal X}$,
yields a
good confidence set $C(X)$ for $\psi$.

However, as pointed out by Casella and Berger (1990) and Casella, Hwang and Robert (1993),
this procedure may lead to very poor confidence sets (confidence sets with
``paradoxical behaviour''). For the remainder of the introduction and in Section 2,
we consider the case that $X = (X_1, \ldots, X_n)$
where $X_1, \ldots, X_n$ are iid $N(\mu, \sigma^2)$ with $\mu$ and $\sigma^2$ unknown,
$\theta = (\mu, \sigma^2)$ and the
parameter of interest is $\mu$. For this case,
Casella, Hwang and Robert (1993) show that, for the standard noninformative prior pdf $\pi(\theta)=1/\sigma^2$
for $\theta$,  the generalized Bayes rule is a very poor confidence interval.
These authors show, however, that the use of the more general class of loss function
\begin{equation*}
\label{new_loss}
L_m(\theta, C) = m \big(\text{length}(C) \big) - {\cal I} (\mu \in C),
\end{equation*}
where $m$ is an appropriately-chosen nonlinear and nondecreasing function, can solve this problem.

In Section 2, we consider the following simple modification of the loss function \eqref{comp_loss_fn}:
\begin{equation}
\label{new_loss}
\tilde{L}(\theta, C) = \frac{\text{length}(C)}{\sigma} - k \, {\cal I}(\mu \in C).
\end{equation}
We show that the standard noninformative prior pdf for $\theta$
leads to a generalized Bayes rule that is the usual confidence interval for $\mu$. In other words,
the ``paradoxical behaviour'' is removed. However, as discussed
in Section 3, we do not advocate the use of generalizations of the loss function \eqref{new_loss} in other contexts.

\bigskip

\noindent {\bf 2. Confidence intervals for the normal mean obtained by using the new loss function \eqref{new_loss}}

\medskip

Suppose that $X_1, \ldots, X_n$ are iid $N(\mu, \sigma^2)$ where both $\mu$ and $\sigma^2$ are unknown ($\mu \in \mathbb{R}$,
$\sigma^2 \in (0, \infty)$).
Let $\theta = (\mu, \sigma^2)$ and suppose that $\mu$ is the parameter of interest. Also let $\bar{X} = \sum_{i=1}^n X_i/n$ and
$S = \sqrt{\sum_{i=1}^n (X_i - \bar{X})^2/(n-1)}$.
Define the quantile $t(m)$ by the requirement that $P \big(-t(m) \le T \le t(m) \big) = 1-\alpha$ for
$T \sim t_m$. The usual $1-\alpha$ confidence interval for $\mu$ is
$\big [ \bar{X} - t(n-1) S/\sqrt{n}, \, \bar{X} + t(n-1) S/\sqrt{n} \big ]$.
Suppose that $\theta$ has the improper prior pdf $\pi(\theta) = 1/\sigma^2$.
This is the standard noninformative prior pdf for $\theta$. Use the new loss function $\tilde{L}(\theta, C)$,
given by \eqref{new_loss}.
In this section, we prove that the generalized Bayes rule is, for the appropriate choice of $k$,
the usual $1-\alpha$ confidence interval for $\mu$.

Since $(\bar{X}, S^2)$ is a sufficient statistic for
$\theta$, we consider confidence intervals for $\mu$ of the form
$C(\bar{X}, S) = \big [ \ell(\bar{X}, S), u(\bar{X}, S) \big]$.
 Define the posterior expected loss
\begin{equation*}
E \big ( \tilde{L}(\theta, C(\bar{X}, S)) \, \big| \, \bar{x}, s \big)
= E \big ( \tilde{L}(\theta, C(\bar{x}, x)) \, \big| \, \bar{x}, s \big),
\end{equation*}
where $E( \, \cdot \, | \, \bar{x}, s)$ denotes the expectation according to the posterior
distribution of $\theta$ i.e. the distribution of $\theta$ conditional on
$(\bar{X}, S) = (\bar{x}, s)$. The posterior expected loss is equal to
\begin{equation}
\label{post_exp_loss_1}
\big( u(\bar{x}, s) - \ell(\bar{x}, s) \big ) \, E(1/\sigma \, | \, \bar{x}, s)
- k \, P \big(\mu \in C(\bar{x}, s) \, | \, \bar{x}, s \big),
\end{equation}
where $P( \, \cdot \, | \, \bar{x}, s)$ denotes the probability according to the posterior
distribution of $\theta$. As is well-known (see e.g. p.215 of Robert, 1994), the marginal
posterior distribution of $\mu$ is such that
\begin{equation*}
\frac{\sqrt{n}(\mu - \bar{x})}{s} \sim t_{n-1}.
\end{equation*}
Thus
\begin{equation*}
P \big(\mu \in C(\bar{x}, s) \, | \, \bar{x}, s \big) =
P \left (  \frac{\sqrt{n}(\ell(\bar{x}, s) - \bar{x})}{s} \le T \le
 \frac{\sqrt{n}(u(\bar{x}, s)- \bar{x})}{s} \right ),
\end{equation*}
where $T \sim t_{n-1}$. As is well-known (see e.g. Box and Tiao, 1973), the marginal
posterior pdf of $\sigma$ is
\begin{equation*}
c(n,s) \sigma^{-n} \exp \left ( - \frac{(n-1)s^2}{2 \sigma^2} \right ),
\end{equation*}
for $\sigma > 0$, where
\begin{equation*}
c(n,s) = \left( \frac{1}{2} \Gamma \left ( \frac{n-1}{2} \right ) \right )^{-1} \left ( \frac{(n-1)s^2}{2} \right )^{(n-1)/2}.
\end{equation*}
Hence
\begin{align*}
E(1/\sigma \, | \, \bar{x}, s) &= c(n,s) \int_0^{\infty} \sigma^{-(n+1)} \exp \left ( - \frac{(n-1)s^2}{2 \sigma^2} \right ) \\
&=c_1(n)/s,
\end{align*}
where
\begin{equation*}
c_1(n) = \frac{\Gamma(n/2)}{\Gamma((n-1)/2)} \sqrt{\frac{2}{n-1}},
\end{equation*}
by (A2.1.4) on p.145 of Box and Tiao (1973). Thus the posterior expected loss \eqref{post_exp_loss_1}
is equal to
\begin{equation*}
\frac{c_1(n) (u(\bar{x},s) - \ell(\bar{x},s))}{s}
- k P \left ( \frac{\sqrt{n}(\ell(\bar{x},s) - \bar{x})}{s} \le T \le \frac{\sqrt{n}(u(\bar{x},s) - \bar{x})}{s} \right).
\end{equation*}
Let $\big(\ell^*(\bar{x},s),u^*(\bar{x},s)\big)$ denote the value of $\big(\ell(\bar{x},s),u(\bar{x},s)\big)$
minimizing the posterior expected loss, subject to $u(\bar{x},s) \ge \ell(\bar{x},s)$. We find this minimizing value as
follows. Define the following function of $(q,r)$:
\begin{equation}
\label{funct_r_q}
\frac{r-q}{s} - k_1(n) P\left( \frac{\sqrt{n} q}{s} \le T \le \frac{\sqrt{n} r}{s} \right ),
\end{equation}
where $k_1(n) = k/c_1(n)$. Let $(q^*, r^*)$ denote the value of $(q,r)$ minimizing \eqref{funct_r_q},
subject to $r \ge q$. Then set $\ell^*(\bar{x},s) = \bar{x} + q^*$ and $u^*(\bar{x},s) = \bar{x} + r^*$.
Let $h = (r-q)/2$ and suppose that $r \ge q$, so that $h \ge 0$. Thus \eqref{funct_r_q} is equal to
\begin{equation}
\label{funct_q_h}
\frac{2h}{s} - k_1(n) P\left( \frac{\sqrt{n} q}{s} \le T \le \frac{\sqrt{n} (q+2h)}{s} \right ).
\end{equation}
We minimize this with respect to $(q,h)$, where $h \ge 0$, in two steps as follows. In the first step,
we minimize \eqref{funct_q_h} with respect to $q$ for fixed $h \ge 0$. We then substitute this minimizing
value of $q$ into \eqref{funct_q_h} and minimize the resulting expression with respect to $h \ge 0$. For
fixed $h \ge 0$, we minimize \eqref{funct_q_h} with respect to $q$ by
maximizing
\begin{equation*}
P \left ( q \le \frac{sT}{\sqrt{n}} \le q+2h \right )
\end{equation*}
with respect to $q$. Clearly, this is maximized by setting $q = -h$. Substituting this value of $q$ into
\eqref{funct_q_h}, we obtain the following function of $h$:
\begin{equation}
\label{funct_h}
\frac{2h}{s} - k_1(n) \left ( 2 F_{n-1} \left( \frac{\sqrt{n} h}{s} \right) - 1 \right),
\end{equation}
where $F_{n-1}$ denotes the $t_{n-1}$ cdf. Multiplying \eqref{funct_h} by $\sqrt{n}$, we obtain
\begin{equation}
\label{funct_h_new}
\frac{2 \sqrt{n} h}{s} - k_2(n) \left ( 2 F_{n-1} \left( \frac{\sqrt{n} h}{s} \right) - 1 \right),
\end{equation}
where $k_2(n) = \sqrt{n} \, k_1(n)$. Minimization of \eqref{funct_h} with respect to $h \ge 0$
is equivalent to minimization of \eqref{funct_h_new} with respect to $h \ge 0$, and this
is equivalent to minimizing
\begin{equation*}
\frac{\sqrt{n} h}{s} - k_2(n)  F_{n-1} \left( \frac{\sqrt{n} h}{s} \right)
\end{equation*}
with respect to $h \ge 0$. Set
\begin{equation*}
k = \frac{1}{f_{n-1}(t(n-1))} \frac{c_1(n)}{\sqrt{n}},
\end{equation*}
where $f_{n-1}$ denotes the $t_{n-1}$ pdf. Thus $k_2(n) = 1/f_{n-1}(t(n-1))$. Our aim, therefore,
is to minimize
\begin{equation*}
g(h) = \frac{\sqrt{n} h}{s} - \frac{1}{f_{n-1}(t(n-1))}  F_{n-1} \left( \frac{\sqrt{n} h}{s} \right)
\end{equation*}
with respect to $h \ge 0$. Now
\begin{align*}
\frac{d g(h)}{dh} &= \frac{\sqrt{n}}{s} - \frac{1}{f_{n-1}(t(n-1))}  f_{n-1} \left( \frac{\sqrt{n} h}{s} \right) \frac{\sqrt{n}}{s} \\
&= \frac{\sqrt{n}}{s} \left(1 - \frac{f_{n-1}(\sqrt{n}h/s)}{f_{n-1}(t(n-1))} \right).
\end{align*}
This derivative is an increasing function of $h \ge 0$ and takes a negative value for $h=0$. Therefore,
$g(h)$ is minimized with respect to $h$ by setting $\sqrt{n}h/s = t(n-1)$, so that
$h = t(n-1)s/\sqrt{n}$ and $C(\bar{x},s) = \big [\bar{x} - t(n-1)s/\sqrt{n}, \bar{x} + t(n-1)s/\sqrt{n} \big]$,
the usual $1-\alpha$ confidence interval for $\mu$.

\bigskip

\noindent {\bf 3. Discussion}

\medskip

The loss function \eqref{new_loss} can be generalized in the obvious way to other contexts
where there is a scaling parameter (analogous to $\sigma$). However, we do not advocate
the use of such a loss function. The expected volume and coverage
probability of a confidence set are very different criteria. An attempt to shoehorn these
criteria into a single risk function that is a linear combination of these criteria
does not seem to be the appropriate strategy.
One is better off to solve the problem of finding a confidence set that minimizes a weighted
average (over the parameter space $\Omega_{\theta}$) of the expected length, subject to the
constraint that this confidence set has coverage probability that never falls below the
specified value $1-\alpha$. In the case that $\theta$ is a scalar and the parameter of
interest, an ingenious solution to this problem is provided by Pratt (1961). Farchione and
Kabaila (2008), Kabaila and Giri (2009ab) solve this problem
in particular settings by computational means.

\bigskip

\noindent {\bf References}

\smallskip

\rf Berger, J.O. 1985. Statistical Decision Theory and Bayesian Analysis, 2nd edition.
Springer-Verlag, New York.

\smallskip

\rf Box, G.E.P., Tiao, G.C. 1973. Bayesian Inference in Statistical Analysis.
Wiley, New York.

\smallskip

\rf   Casella, G., Berger, R.L., 1990. Statistical Inference. Wadsworth, Belmont, CA.

\smallskip

\rf   Casella, G., Hwang, J.T.G., Robert, C., 1993. A paradox in decision-theoretic
interval estimation. Statistica Sinica, 3, 141--155.

\smallskip

\rf   Farchione, D., Kabaila, P., 2008. Confidence intervals for the
normal mean utilizing prior information. Statistics \& Probability Letters
78, 1094--1100.

\smallskip

\rf Kabaila, P., Giri, K., 2009a. Confidence intervals in regression
utilizing uncertain prior information. Journal of Statistical Planning and
Inference 139, 3419--3429.

\newpage


\rf Kabaila, P., Giri, K., 2009b. Large-sample confidence intervals for the treatment difference in a two-period
crossover trial, utilizing prior information. Statistics \& Probability Letters 79, 652--658.

\smallskip

\rf Lehmann, E.L., Casella, G. 1998. Theory of Point Estimation, 2nd edition.
Springer-Verlag, New York.

\smallskip

\rf Robert, C.P. 1994. The Bayesian Choice, A Decision-Theoretic Motivation.
Springer-Verlag, New York.

\smallskip

\rf Pratt, J.W. 1961. Length of confidence intervals. Journal of the American Statistical Association
56, 549--567.

\end{document}